\documentclass{birkjour}

\usepackage[T2A,T1]{fontenc}
\usepackage[utf8]{inputenc}
\usepackage[russian,english]{babel}
\usepackage{epigraph}
\usepackage{skull}

\title[Remembering Victor Petrovich Havin]
{Remembering Victor Petrovich Havin}

\author[K. M. Dyakonov]{Konstantin M. Dyakonov}

\address
{Departament de Matem\`atiques i Inform\`atica\\ 
Universitat de Barcelona, IMUB and BGSMath\\
Gran Via, 585\\ 
E-08007 Barcelona\\ 
Spain}

\address{\quad\\
Instituci\'o Catalana de Recerca i Estudis Avan\c{c}ats (ICREA)\\ 
Pg. Llu\'is Companys, 23\\ 
E-08010 Barcelona\\ 
Spain}

\email{konstantin.dyakonov@icrea.cat}
\thanks{Supported in part by grant MTM2014-51834-P from El Ministerio de 
Econom\'ia y Competitividad (Spain) and grant 2014-SGR-289 from AGAUR (Generalitat de Catalunya).}

\begin{document}

\maketitle

\epigraph{\foreignlanguage{russian}{Когда погребают эпоху,\\
Надгробный псалом не звучит.}}
{\foreignlanguage{russian}{Анна Ахматова}}

Victor Petrovich Havin will probably be remembered for many things. Of course, he was an eminent mathematician -- quite a towering figure in his field -- and the founding father of a renowned analysis school. He was also a wonderful lecturer (certainly the best I've ever listened to), a truly charismatic teacher, and a man of great personal charm. This last quality is, in fact, crucial in understanding why he was able to attract and influence so many; yet it seems to be the least describable part, too. Indeed, that mysterious (though unmistakably recognizable) substance known as personal charm -- or glamour, or charisma, whatever -- defies description. Why do some people have it, while others have none? God knows, but he won't tell. 

\par I find it equally hard to explain why V.P.'s lectures were so inspiring, but I have a clear memory of being immediately fascinated by his personality and his subject, Mathematical Analysis, when I first heard him speak. (That was 35 years ago, so he was 48, younger than I am now. The audience consisted of over 100 freshers aged 17 plus epsilon, crammed in a large lecture theatre at {\it Mat-Mekh}, in Old Peterhof.) Somehow, he taught {\it as one who has authority}, and you had the feeling that if this highly intelligent, eloquent and overall extraordinary person is so excited about the matter at hand, then that's probably the sort of stuff you should be doing -- what's good enough for Havin is good enough for me. 

\par The mathematical content of V.P.'s lectures would often be quite technical: a host of $\varepsilon$'s and $\delta$'s at the initial stage, then lots of {\it hard analysis} proofs involving lengthy calculations, combinatorial arguments, sophisticated estimates of various unwieldy integrals and whatnot at a more advanced level. (Quite frankly, at times his heavy notation and terminology seemed somewhat too dense to me, and I found myself thinking that certain terms -- such as {\it $\varepsilon$-\foreignlanguage{russian}{допуск}}, say -- could be safely dispensed with.) However, even those technical matters were presented in a truly artistic manner, everything was extremely well organized, and you never felt bored. On the contrary, it was always clear that something interesting and important was going on. 

\par Besides, whenever possible, V.P. would begin with an informal and/or heuristic \lq\lq hand-waving argument\rq\rq\, to unveil the simple truth ({\it \foreignlanguage{russian}{сермяжная правда}}, he used to say) underlying this or that result; the formal proof that ensued would then become much more digestible. For example, to convince us that $f\circ g$ is continuous provided both $f$ and $g$ are, V.P. would first write down something like 
$$a\approx b\implies g(a)\approx g(b)\implies f(g(a))\approx f(g(b))$$
by way of explanation, and then convert this \lq\lq quasi-proof\rq\rq\, into
a rigorous $\varepsilon$-$\delta$ argument. Likewise, brief informal comments were offered to clarify the implicit function theorem (whose full proof was quite long and exhausting, if I recall well), the change of variable theorem (involving the Jacobian) for integrals in $\mathbb R^n$, and other real analysis facts. Remarkably, though, when it came to Complex Analysis, V.P. confessed that he could see no simple explanation for Cauchy's theorem (the one that says
$$\oint_{\partial\Omega}f(z)dz=0$$
for $f$ holomorphic on $\overline\Omega$), no {\it homespun truth} behind it that would make the result essentially obvious. Rather, he viewed it as an amazing fact with an air of magic about it, a miracle that continued to surprise him ever since he had come across it a few decades before. Now he was genuinely elated at rediscovering this and other gems of the complex realm (for that's where his heart belonged), and he knew how to convey elation to his disciples. 

\par At the same time, it wasn't Havin's idea to conquer the audience by making the presentation entertaining. After all, he was there to provide us with food for thought, and the student was supposed to toil in the sweat of his/her brow (this last point was repeatedly emphasized). And not only that, the student was assumed to be obsessed with mathematics above all other things! Needless to say, this was not always the case in reality, and occasionally the students' ignorance and/or apathy seemed to grind V.P.'s gears. He would then address some bitter words of disappointment to those guilty (or, by extension, to all those he had in front). 

\par As any lecturer, V.P. had some home-made jokes in his toolbox, but not many. For instance, he would warn us -- with a straight face -- that \lq\lq the triangle inequality\rq\rq\, ({\it\foreignlanguage{russian}{неравенство треугольника}}) should not be misinterpreted as \lq\lq Triangle's inequality,\rq\rq\, i.e., a result due to the hypothetical mathematician named Triangle ({\it\foreignlanguage{russian}{Треугольник}}), rather an unlikely interpretation indeed. Many an important formula on the blackboard was boxed, supplied with a skull and crossbones symbol $\big(\skull\big)$, and much woe was prophesied to those poor buggers among us who would even begin to hesitate, albeit for a tiny fraction of a second, in reproducing it at the exam if required. Any such student would be mercilessly flunked, V.P. told us in a menacing tone. The formulas you had to forget in order to deserve such capital punishment included the basic Taylor--Maclaurin expansions of elementary functions, the Euler identities 
$$\cos z=\frac{e^{iz}+e^{-iz}}2,\qquad\sin z=\frac{e^{iz}-e^{-iz}}{2i}$$ 
and so forth. 

\par What really mattered, though, was V.P.'s unique emotional style, his individualistic -- sometimes strikingly original -- view on many mathematical phenomena, and the highly expressive metaphorical language he used to describe them. For instance, the zeros of analytic functions were said to be {\it contagious} (if there are too many, they propagate like plague until they absorb everything), and the B\"urmann--Lagrange inversion formula was claimed to make the {\it golden dream of mankind} come true (namely, that every equation be solvable). Conversely, V.P. would occasionally resort to mathematical language to explain certain phenomena outside mathematics. One such digression, related again to the identity principle for analytic functions, concerned the nature of great works of art (e.g., in poetry or music), and V.P.'s point was that {\it analyticity} is a key feature of any true masterpiece. Indeed, once you have a germ (say, the first line of a poem or a few opening bars of a musical piece), it is bound to extend on its own, and in the only possible way; the continuation is unique, and you have no further control over it. Thus, according to V.P., if it begins with 

\medskip
\begin{center}
{\it\foreignlanguage{russian}{Буря мглою небо кроет,}}
\end{center}

\medskip\noindent
then the next line, and in fact the rest of the poem, is predetermined: it just {\it has to be} what it is. Consequently, the poet doesn't have to sit there thinking what to write next, for he is but a tool employed by that divine force, the analytic extension, to do the job. 

\par Havin's own language was poetic -- in his lectures, writings, and daily life -- and pertinent quotations from (or allusions to) various authors would often serve him to reinforce it. For instance, when running out of time and forced to skip certain details, he said he was \lq\lq setting his heel on the throat of his own song,\rq\rq\, an obvious reference to Mayakovsky's poem \lq\lq At the Top of My Voice\rq\rq. On another occasion, and unexpectedly enough, V.P. borrowed Mayakovsky's lines 

\medskip
\begin{center}
{\it\foreignlanguage{russian}{Но поэзия -- \qquad\qquad\qquad\qquad\\
\qquad\qquad\qquad\qquad\qquad\qquad пресволочнейшая штуковина: \\
существует -- \qquad\qquad\qquad\qquad\\
\qquad\qquad\qquad\qquad и ни в зуб ногой}}
\end{center}

\medskip\noindent
to describe the evasive nature of mathematics (the only change suggested was to replace {\it\foreignlanguage{russian}{поэзия}} by {\it\foreignlanguage{russian}{математика}}); actually, this was incorporated into a toast that V.P. made at a banquet. 

\par Elsewhere -- namely, in one of his appendices to Paul Koosis' book \lq\lq Introduction to $H_p$ Spaces\rq\rq\, -- V.P. invoked Tyutchev's famous line

\medskip
\begin{center}
{\it\foreignlanguage{russian}{Мысль изреченная есть ложь}}
\end{center}

\medskip\noindent
(\lq\lq A thought, once uttered, is a lie\rq\rq) to support his feeling that the essence of a theorem can seldom, if ever, be grasped from its statement alone. This remarkable claim seems to be typical of Havin's mathematical philosophy, and I take the liberty of reproducing a bit of his text in Koosis' translation: 
{\it\lq\lq Confined within the narrow limits of its formulation, a theorem does not tell {\rm all} of the truth about itself and is perhaps only fit for inclusion in a handbook. Its true meaning is inseparable from the proof (or proofs, if there are several of them).\rq\rq}

\par One quotation that comes to mind in connection with Havin himself is Weierstrass' often cited maxim that \lq\lq a mathematician who is not somewhat of a poet, will never be a perfect mathematician\rq\rq. Well, V.P. undoubtedly met the criterion. 

\par It was my good fortune to be a Ph.D. student of Havin in my post-university years. Strictly speaking, I was rather a {\it\foreignlanguage{russian}{соискатель}}, a fairly vague status implying that I had no formal link to the university, and hence to V.P., while working on my thesis. Nonetheless, he took care of me with great kindness, and his fatherly attitude went far beyond what a supervisor is normally expected to offer. Once in a while I would ring him up, then come to his place, we would sit down at his enormous desk with piles of books, manuscripts, reprints, preprints, letters, photos and whatnot on it, and talk about mathematics. This could go on for hours. I would tell him about any progress I was (or wasn't) able to make in the past weeks, and his vivid reaction -- sometimes approving, even enthusiastic, sometimes perhaps tinged with skepticism -- meant a lot to me. V.P.'s excellent taste and keen intuition were as impressive as his great fund of knowledge, and talking to him would broaden your horizons. Moreover, he had that special flair, a kind of perfect pitch in mathematics, which allowed him to tell you, with a reliable degree of certainty, whether your stuff was attractive and likely to have interesting connections with other things. Some of the insights and associations that occurred to him along the way were really illuminating. 

\par As a \lq\lq mathematical son\rq\rq\, of V.P., you could always count on his advice, reassurance (in case you got stuck) and -- most importantly -- encouragement. Yes, encouragement is what you need badly when you are a beginner with little or no self-confidence, with no achievements yet to boast of; and V.P. did his best to encourage you generously whenever you happened to make a perceptible step forward. Sometimes you were lucky enough to do more -- crack the problem completely, or perhaps chance upon an unexpected twist that led to a sexy result -- and V.P. would be really pleased. He would speak highly of your result, then probably do some advertising by telling other people about it and/or getting you to speak at the seminar... To see him excited or surprised by what you had done was very gratifying indeed, and his praise made you feel euphoric. 

\par I did not see much of V.P. after I had moved to Spain (which happened around 2000), but I kept paying him a visit once a year, usually towards the end of August, while on vacation in St. Petersburg. V.P. and Valentina Afanas'evna, his wife, would kindly offer me a cup of coffee or a glass of wine; then V.P. and I would spend a few hours talking in his room (following the tradition, and perhaps abusing his friendly attention, I went on testing my new theorems on him); alternatively, weather permitting, we might go for a walk. As always, V.P. irradiated that peculiar warmth which made his company so enjoyable; it was invariably a pleasure to be by his side and listen to what he had to say. When asked about his health -- which left very much to be desired -- V.P. would respond briefly and switch to some other topic. He would speak of the current state of things at {\it Mat-Mekh} (being often displeased with the changes that modern times had brought about), of what it was like in former days (his reminiscences and anecdotes involved great names of the past, such as V.I. Smirnov, V.A. Rokhlin, A.D. Alexandrov, S.N. Mergelyan), or he might touch upon his own recent work, usually done in collaboration with this or that Ph.D. student. During my last visit, in August 2015, he explained to me -- in a fairly detailed way -- his joint paper with Pavel Mozolyako on the normal variation of a positive harmonic function. 

\par On a less serious note, I recall a quasi- (or pseudo?) mathematical puzzle, which V.P. once offered me by way of entertainment: compute the product 
$$\sin\alpha\cdot\sin\beta\cdot\sin\gamma\cdot...\cdot\sin\omega.$$ 
(Spoiler: the answer is $0$, since one of the factors is $\sin\pi$.)

\par Further topics of our conversations -- or his monologues -- included politics, literature, and most notably, music. In fact, V.P. had a life-long passion for music and spoke about it knowledgeably. He admired Mahler's symphonies and Shostakovich's quartets, knew a large number of operas thoroughly, if not by heart. I also remember him saying that Beethoven's even-numbered symphonies deserve more attention than the (overwhelmingly more famous) odd-numbered ones; V.P. traced this opinion back to V.I. Smirnov, his own mathematical father, and seemed to concur with it. I, for one, still don't feel fully convinced... 

\par By the way, V.P. had a fine bass voice and sang beautifully. His repertoire seemed to consist chiefly of {\it\foreignlanguage{russian}{блатные песни}} -- the songs attributed to, and describing the miserable life of, the criminal dregs of society. (I have always suspected that this kind of folk culture is largely produced by intellectuals, maybe even distinguished professors.) Although I had no chance to witness any live performance by V.P., I did have an occasion to listen to him on tape. The record was provided by the late Evsey Matveevich (Seva) Dyn'kin when I was visiting him in Haifa, back in 1998, and he played it on his tape recorder for me. I had a strange feeling when V.P.'s voice entered the room: the genre and subject matter were somewhat unusual (I would have rather expected V.P. to besing harmonic vector fields or the like), but after all, his account of a bank holdup was no less emotional: 

\medskip
\begin{center}
{\it\foreignlanguage{russian}{\qquad\quad Взяли в сберкассе мы сумму немалую -- \\
Двадцать пять тысяч рублей...}}
\end{center}

\medskip
\par He was a lovable and unforgettable human being. We shall not look upon his like again.

\end{document}